%% file: sinkhorn_geometry.tex
\documentclass[a4paper,10pt,reqno]{amsart}

\usepackage[utf8]{inputenc}
\usepackage[
theoremdefs,
final
]{latexdev}

\usepackage[numbers]{natbib}
\usepackage{standalone}
\usepackage{tikz,tikz-cd}

\usepackage{algpseudocode}
\usepackage{algorithm}

\input{definitions}

\title[Geometry of the Sinkhorn algorithm]{On the geometry and dynamical formulation of the Sinkhorn algorithm for optimal transport}
\author[Klas Modin]{Klas Modin}
\address{Department of Mathematical Sciences, Chalmers University of Technology and University of Gothenburg, SE-412 96 Gothenburg, Sweden}
\email{klas.modin@chalmers.se}
\date{\today}                                           

\dedicatory{On the occasion of the 60th birthdays of Hans Munthe--Kaas and Brynjulf Owren}

\begin{document}

\begin{abstract}
The Sinkhorn algorithm is a numerical method for the solution of optimal transport problems. Here, I give a brief survey of this algorithm, with a strong emphasis on its geometric origin: it is natural to view it as a discretization, by standard methods, of a non-linear integral equation. In the appendix, I also provide a short summary of an early result of Beurling on product measures, directly related to the Sinkhorn algorithm.


\end{abstract}

\maketitle


\section{Introduction}\label{sec:intro}

The Sinkhorn algorithm has quickly sailed up as a popular numerical method for optimal transport (OT) problems~(see the book by Peyré and Cuturi~\cite{PeCu2020} and references therein).
Its standard derivation starts from the Kantorovich formulation of the discrete OT problem and then adds entropic regularization, which enables Sinkhorn's theorem~\cite{Si1964} for the optimal coupling matrix (see the paper by Cuturi~\cite{Cu2013} for basic notions, and the work of Schmitzer~\cite{Sc2019} and the PhD thesis of Feydy~\cite{Fe2020_thesis} for more details, including optimized and stable implementation of the full algorithm).
In this paper, I give a brief survey of the Sinkhorn algorithm from a geometric viewpoint, making explicit connections to well-known techniques in geometric hydrodynamics, such as the Madelung transform, and the dynamical formulation of OT as given by Benamou and Brenier~\cite{BeBr2000}.
This viewpoint on entropic regularization is known to experts since long before the Sinkhorn algorithm became popular for OT: it was first formulated by Schrödinger~\cite{Sc1931}, then re-established by Zambrini~\cite{Za1986}, and is today well-known as the law governing the \emph{Schrödinger bridge problem} in probability theory (cf.~Leonard~\cite{Le2014}).
Here, I wish to complement this presentation by entirely avoiding probability theory and instead use the language of geometric analysis, particularly geometric hydrodynamics (\emph{cf}.\ Arnold and Khesin~\cite{ArKh1998}).
The presentation of Léger~\cite{Le2019}, based on Wasserstein geometry, is close to mine and contains more details.
A benefit of the hydrodynamic formulation is that it readily enables standard numerical tools for space and time discretization (and the accompanying numerical analysis).

The principal theme of my presentation is that the Sinkhorn algorithm on a smooth, connected, orientable Riemannian manifold $M$ has a natural interpretation as a space and time discretization of the following non-linear evolution integral equation:
\begin{equation}\label{eq:integral_equations}
	\left\{
	\begin{aligned}
		\partial_s f &= - f -  \log\left(\int_M K_\varepsilon(\cdot, y)\, \ee^{g(y)}\uud y \right) + \log\rho_0 ,\qquad f\colon M\times[0,\infty)\to \RR \\
		\partial_s g &= - g -  \log\left(\int_M K_\varepsilon(\cdot, y)\, \ee^{f(y)}\uud y \right) + \log\rho_1 ,\qquad g\colon M\times[0,\infty)\to \RR .
	\end{aligned}
	\right.
\end{equation}
Here, $\rho_0$ and $\rho_1$ are two strictly positive and smooth densities with the same total mass relative to the Riemannian volume form, and $K_\varepsilon(x,y)$ is the heat kernel on $M$.
In the special case $M=\RR^n$ then
\begin{equation}\label{eq:heat_kernel_Rn}
	K_\varepsilon(x,y) = \frac{1}{(4\pi \varepsilon)^{n/2}}\exp\Big(-\frac{\norm{x-y}^2}{4\varepsilon}\Big) .
\end{equation}
It is sometimes convenient to use the heat flow semi-group notation
\begin{equation}
	\ee^{\varepsilon\Delta}f \coloneqq \int_M K_\varepsilon(\cdot,y)f(y)\ud y.
\end{equation}

\begin{remark}
	The equations \eqref{eq:integral_equations} are almost linear for small $\varepsilon$.
	Indeed, they can be written
	\begin{equation}
	\left\{
	\begin{aligned}
		\partial_s f &= - f - g + \log\rho_0 - T_\varepsilon(g), \\
		\partial_s g &= - g - f + \log\rho_1 - T_\varepsilon(f),
	\end{aligned}
	\right.
	\qquad T_\varepsilon(f) \coloneqq \log \big(\ee^{-f}\ee^{\varepsilon\Delta} \ee^{f}\big), 
\end{equation}
where $(f,\varepsilon) \mapsto T_\varepsilon(f)$ is a $C^\infty$ mapping, for example in Sobolev topologies $H^s$. 
Local existence and uniqueness of solutions to the integral equations~\eqref{eq:integral_equations} are obtained by standard Picard iterations.
The extension to global results probably follow by standard techniques, but the question of \emph{convergence} as $s\to\infty$ (and possibly $\epsilon\to 0$ simultaneously) is subtle.
To approach it, one could attempt backward error analysis on the continuous version of the system \eqref{eq:splitting_on_M_fg} below.
I expect this would yield a connection to the non-linear parabolic equation studied by Berman~\cite{Be2020}, whose dynamics, he proved, approximates the dynamics of $\vect{f}^{(k)}\mapsto \vect{f}^{(k+1)}$ below.
\end{remark}

Before discussing the geometric origin of the equations \eqref{eq:integral_equations}, let us see how the standard (discrete) Sinkhorn algorithm arise from the equations \eqref{eq:integral_equations}.

Let $\rho_0$ and $\rho_1$ be two atomic measures on $M$
\begin{equation}\label{eq:discrete_marginals}
	\rho_0 = \sum_{i=1}^N p_i \delta_{x_i}, \qquad \rho_1 = \sum_{i=1}^N q_i \delta_{y_i},
\end{equation}
where $x_i,y_i \in \RR^n$ and the weights $p_i>0$ and $q_i>0$ fulfill $\sum_i p_i = \sum_i q_i$.
If we change variables $a(x, s) =  \exp(f(x,s))$ and $b(x,s) = \exp(g(x,s))$ then the integral equations~\eqref{eq:integral_equations} become
\begin{equation}\label{eq:integral_equations_ab}
	\left\{
	\begin{aligned}
		\partial_s a &= -a \log\left(\frac{a\ee^{\varepsilon \Delta}b}{\rho_0} \right), \\
		\partial_s b &= -b \log\left(\frac{b\ee^{\varepsilon \Delta}a}{\rho_1} \right)  .
	\end{aligned}
	\right.
\end{equation}
For this to make sense for~\eqref{eq:discrete_marginals} we need $a \ll \rho_0$ and $b \ll \rho_1$, i.e., $a=\sum_i a_i \delta_{x_i}$ and $b=\sum_i b_i\delta_{y_i}$ for some $a_i \geq 0$ and $b_i \geq 0$.
Since the heat flow convolves a delta distribution to a smooth function, it follows if $\varepsilon>0$ that
\begin{equation}
	\left(\frac{a\ee^{\varepsilon \Delta}b}{\rho_0}\right)(x_i) = \frac{a_i}{q_i}\sum_{j=1}^N K_\varepsilon(x_i, y_j) b_j = \frac{a_i}{q_i} \vect{K}_\varepsilon\vect{b},
\end{equation}
where $\vect{b} = (b_1,\ldots,b_N)$ and $\vect{K}_\varepsilon$ is the $N\times N$ matrix with entries $(\vect{K}_\varepsilon)_{ij} = K_\varepsilon(x_i,y_i)$.
Likewise,
\begin{equation}
	\left(\frac{b\ee^{\varepsilon \Delta}a}{\rho_1}\right)(y_i) = \frac{b_i}{p_i} \vect{K}_\varepsilon^\top\vect{a},
\end{equation}
where $\vect{a} = (a_1,\ldots,a_N)$ and the heat kernel property $K_\varepsilon(y_i, x_j) = K_\varepsilon(x_j, y_i)$ is used.
Expressed in $\vect{a}$ and $\vect{b}$ the equations \eqref{eq:integral_equations_ab} now become an ordinary differential equation (ODE) on $M^{2N}$
\begin{equation}\label{eq:ode_on_M}
	\left\{
	\begin{aligned}
		\partial_s{\vect{a}} &= -\vect{a}* \log\left( \frac{\vect{a}}{\vect{q}}* \vect{K}_\varepsilon\vect{b} \right), \\
		\partial_s{\vect{b}} &= -\vect{b}* \log\left( \frac{\vect{b}}{\vect{p}}* \vect{K}_\varepsilon^\top\vect{a} \right),
	\end{aligned}
	\right.
\end{equation}
where $*$ denotes element-wise multiplication and $\log$ and divisions are also applied element-wise.
Notice that this ODE loses its meaning (or rather its connection to \eqref{eq:integral_equations_ab}) if $\varepsilon = 0$.
Moving back to the coordinates $\vect{f} = \log\vect{a}$ and $\vect{g} = \log\vect{b}$ yields the system
\begin{equation}\label{eq:ode_on_M_fg}
	\left\{
	\begin{aligned}
		\partial_s{\vect{f}} &= -\vect{f} - \log\left(\vect{K}_\varepsilon\exp\vect{g} \right) + \log\vect{q}, \\
		\partial_s{\vect{g}} &= -\vect{g} - \log\left(\vect{K}_\varepsilon^\top\exp\vect{f} \right) + \log\vect{p},
	\end{aligned}
	\right.
\end{equation}

To proceed, we need a time discretization. 
For this, apply to \eqref{eq:ode_on_M_fg} the Trotter splitting (\emph{cf}.~\cite{McQu2002}) combined with the forward Euler method to obtain
\begin{equation}\label{eq:splitting_on_M_fg}
	\left\{
	\begin{aligned}
		\vect{f}^{(k+1)} &= (1-h)\vect{f}^{(k)} - h\log\left(\vect{K}_\varepsilon\exp\vect{g}^{(k)} \right) + h\log\vect{q}, \\
		\vect{g}^{(k+1)} &= (1-h)\vect{g}^{(k)} - h\log\left(\vect{K}_\varepsilon^\top\exp\vect{f}^{(k+1)} \right) + h\log\vect{p},
	\end{aligned}
	\right.
\end{equation}
where $h>0$ is the time-step length.

\begin{definition}\label{def:sinkhorn_on_M}
	The \emph{discrete Sinkhorn algorithm on $M$} is given by the time discretization~\eqref{eq:splitting_on_M_fg} with $h=1$.
\end{definition}

\begin{remark}
	For $M=\RR^n$ the heat kernel is given by \eqref{eq:heat_kernel_Rn}.
	If we express \eqref{eq:splitting_on_M_fg} in the variables $\vect{a}$ and $\vect{b}$ and take $h=1$ we recover the Euclidean Sinkhorn algorithm as presented in the literature
	\begin{equation}\label{eq:standard_sinkhorn}
		\vect{a}^{(k+1)} = \frac{\vect{q}}{\vect{K}_\varepsilon\vect{b}^{(k)}}, \qquad \vect{b}^{(k+1)} = \frac{\vect{p}}{\vect{K}_\varepsilon^\top\vect{a}^{(k+1)}}.
	\end{equation}
\end{remark}

\begin{figure}
	\raggedright
	\includegraphics[width=0.9\textwidth]{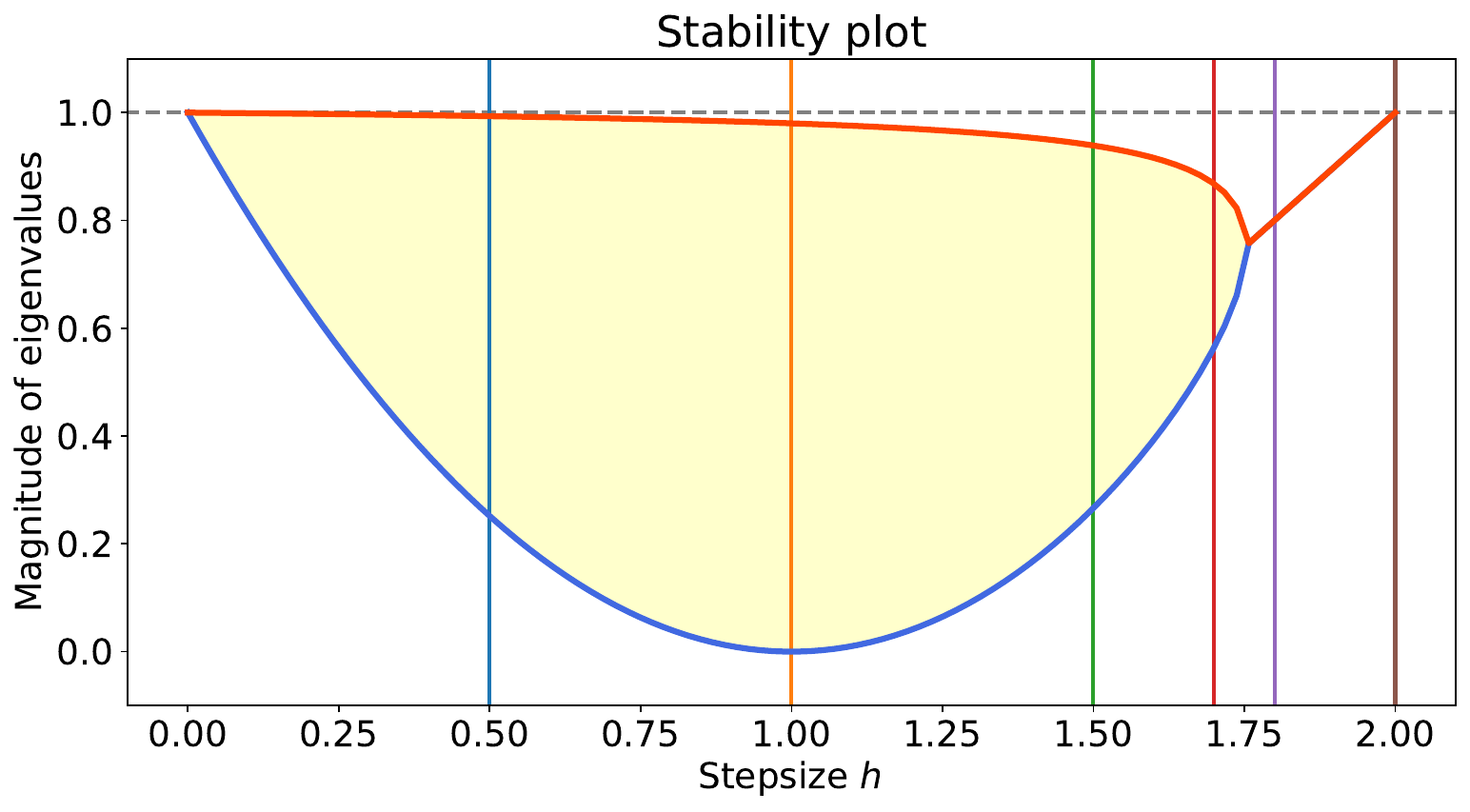}
	\includegraphics[width=0.9\textwidth]{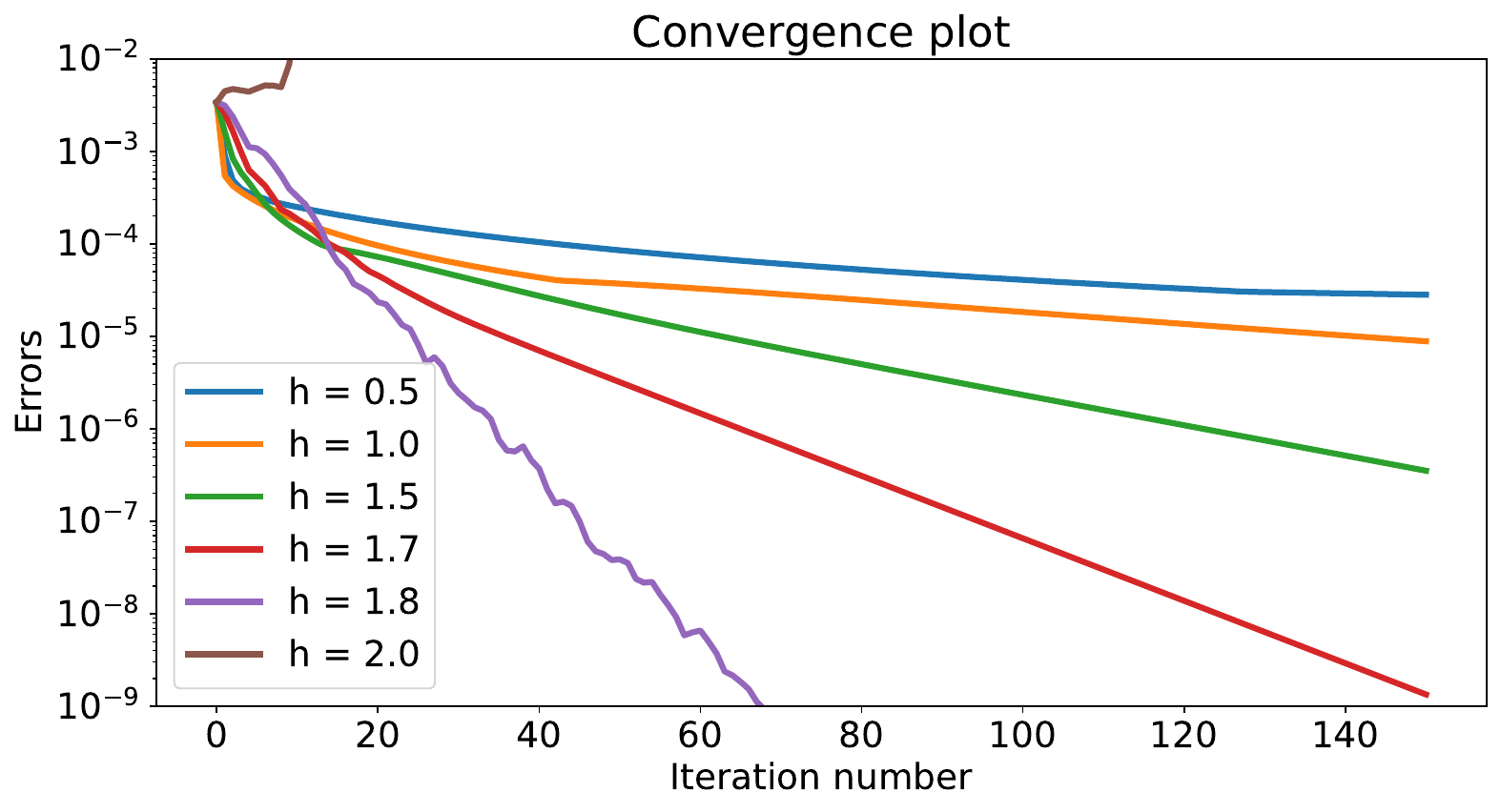}
	\caption{
		Stability analysis (upper figure) for the Trotter--Euler splitting method applied to the linear test equation $\dot x = -x - (1-\delta)y, \; \dot y = -(1-\delta)x - y$ for $\delta = 10^{-2}$, roughly corresponding to the system \eqref{eq:splitting_on_M_fg} with $\varepsilon \simeq \delta$.
		The method is stable whenever the magnitude of the corresponding two eigenvalues for the flow map are bounded by 1. 
		The lower figure shows, for $\varepsilon = 10^{-2}$ and various step-sizes marked in the stability plot, how the $L^2$ norm of the right-hand side of \eqref{eq:ode_on_M_fg} (the error) decreases with the number of time-steps taken with the Trotter--Euler method \eqref{eq:splitting_on_M_fg}.
		The convergence plot matches the theoretical stability plot well.
		At $h=1$ we observe an initially almost vertical decrease of the error due to the small eigenvalue near zero, thereafter the dynamics takes place in the eigenspace of the other eigenvalue close to one and here the convergence is slow.
		The optimal step-size is where the two eigenvalues meet at about $h\approx 1.75$.
		In the interval $h\in [1.75, 2)$ we expect an almost linear decrease of the error, as the two eigenvalues are the same. 
		For $h\geq 2$ the method is not (linearly) stable.
		}\label{fig:stability}
\end{figure}

\begin{remark}
	If we take $h\neq 1$ then \eqref{eq:splitting_on_M_fg} expressed in $\vect{a}$ and $\vect{b}$ become
	\begin{equation}\label{eq:overrelaxed_sinkhorn}
		\vect{a}^{(k+1)} = \left(\vect{a}^{(k)} \right)^{1-h}\left(\frac{\vect{q}}{\vect{K}_\varepsilon\vect{b}^{(k)}}\right)^h, \qquad \vect{b}^{(k+1)} = \left(\vect{b}^{(k)} \right)^{1-h}\left(\frac{\vect{p}}{\vect{K}_\varepsilon^\top\vect{a}^{(k+1)}}\right)^h.
	\end{equation}
	For $h>1$ this is the \emph{over-relaxed Sinkhorn algorithm}, which convergences faster than~\eqref{eq:standard_sinkhorn} (see~\cite{ThChDoPa2021,PeChViSo2019,LeReSaUs2021}).
	Indeed, the faster convergence is readily understood by applying linear stability theory to \eqref{eq:splitting_on_M_fg} in the vicinity $\varepsilon \to 0^+$. 
	From a numerical ODE perspective, the accelerated convergence is expected: larger time steps typically yield faster marching toward asymptotics, provided that the time step is small enough for the method to remain stable.
	Indeed, Figure~\ref{fig:stability} shows an almost perfect match between the stability region of the Trotter-Euler method and the convergence of the time-step iterations~\eqref{eq:splitting_on_M_fg}.
	An interesting venue could be to look for other methods, with other numerical damping properties.
\end{remark}

\begin{remark}
	Notice that if $C = \int_M a\ee^{\varepsilon\Delta}b$ then, along the flow~\eqref{eq:integral_equations_ab},
	\begin{equation*}
		\frac{d}{ds} C = - \int_M a \ee^{\varepsilon\Delta}b \log\left(\frac{a \ee^{\varepsilon\Delta}b}{\rho_0} \right)
		- \int_M b \ee^{\varepsilon\Delta}a \log\left(\frac{b \ee^{\varepsilon\Delta}a}{\rho_1} \right).
	\end{equation*}
\end{remark}

\subsection*{Acknowledgement}
This work was supported by the Swedish Research Council (grant number 2022-03453) and the Knut and Alice Wallenberg Foundation (grant number WAF2019.0201).
I would like to thank Ana Bela Cruzeiro, Christian Leonard, and Jean-Claude Zambrini, for helpful and intriguing discussions, and especially for pointing me to the ``forgotten'' work of Beurling.  

\section{Dynamical formulation of optimal transport}\label{sec:bridge}

This section reviews the dynamical (or fluid) formulation of smooth optimal transport problems as advocated by Benamou and Brenier~\cite{BeBr2000}.
See~\cite{KhMiMo2019} for details on the notation and more information about infinite-dimensional manifold and Riemannian structures.
For simplicity, I assume from here on that $M$ is a compact Riemannian manifold without boundary.
The non-compact or boundary cases can be handled by introducing suitable decay or boundary conditions.

Let $\Dens(M) = \{ \rho \in C^\infty(M) \mid \rho(x) > 0,\, \int_M \rho = 1\}$ denote the space of smooth probability densities.
It has the structure of a smooth Fréchet manifold~\cite{Ha1982}.
Its tangent bundle is given by tuples $(\rho,\dot\rho)$ where $\dot\rho \in C^\infty_0(M) = \{ f \in C^\infty(M)\mid \int_M f = 0\}$.
Otto~\cite{Ot2001} suggested the following (weak) Riemannian metric on $\Dens(M)$
\begin{equation}\label{eq:otto_metric}
	\pair{\dot\rho,\dot\rho}_{\rho} = \int_M \abs{\nabla S}^2\rho, \qquad \dot\rho + \divv(\rho\nabla S) = 0.
\end{equation}
The beauty of this metric is that the distance it induces is exactly the $L^2$-Wasserstein distance, and the geodesic two-point boundary value problem corresponds to the optimal transport problem (in the smooth category).
See \cite{Ot2001,KhWe2009,Mo2017} for details.
In summary, $L^2$ optimal transport is a problem of Lagrangian (variational) mechanics on $\Dens(M)$: find a path $[0,1] \ni t \mapsto \rho_t\in\Dens(M)$ with fixed end-points $\rho_0$ and $\rho_1$ that extremizes (in this case minimizes) the action functional
\begin{equation}\label{eq:action}
	A(\rho_t) = \int_{0}^1 L\Big(\rho_t,\frac{\partial \rho_t}{\partial t}\Big)\, \ud t,
\end{equation}
for the kinetic energy Lagrangian $L(\rho,\dot\rho) = \frac{1}{2}\pair{\dot\rho,\dot\rho}_{\rho}$.
The optimal transport map is then recovered as the time-one flow map for the time dependent vector field on $M$ given by $v(x,t) = \nabla S_t(x)$.

The dynamical formulation is now obtain by replacing the variational problem for the action \eqref{eq:action} with an equivalent constrained variational problem on $\Dens(M)\times \Omega^1(M)$ (densities times one-forms) for the action
\begin{equation}
	\bar A(\rho_t,m_t) = \frac{1}{2}\int_0^{1} \inner{m_t/\rho_t, m_t}_{L^2}\, \ud t
\end{equation}
subject to the constraint
\begin{equation}
	\dot\rho_t + \divv(m_t) = 0 .
\end{equation}
This is a convex optimization problem since $\bar A$ is convex and the constraint is linear (see \cite{BeBr2000} for details).
Notice that the convexity of $\bar A$ is with respect to the linear structure of $C^\infty(M)\times \Omega^1(M)$, which is different from the non-linear convexity notion for the Levi-Civita connection associated with the Riemannian metric \eqref{eq:otto_metric}.

\section{Entropic regularization and the Madelung transform}\label{sec:hopfcole}

The aim of this section is to introduce entropic regularization of the dynamical formulation and show how it simplifies the problem via the imaginary Madelung (or Hopf--Cole) transform.
This transformation is the analog, in the dynamical formulation of smooth OT, to Sinkhorn's theorem applied to the coupling matrix in discrete OT.
Most of the results presented in this section are available in the papers by Leonhard~\cite{Le2014} and by Léger~\cite{Le2019}.
See also Léger and Li~\cite{LeLi2021} for a generalized Schrödinger bridge problem.

Let me first introduce two central functionals from information theory. 
The first is \emph{entropy}, i.e., the functional on $\Dens(M)$ given by
\begin{equation}
	E(\rho) = \int_M \rho \log\rho .
\end{equation}
Its cousin, the \emph{Fisher information functional}, is given by
\begin{equation}
	I(\rho) = \frac{1}{2}\int_M \frac{\abs{\nabla \rho}^2}{\rho}.
\end{equation}
There are various ways to describe the relation between $E(\rho)$ and $I(\rho)$:
\begin{itemize}
	\item $I(\rho)$ is the trace of the Hessian (with respect to~\eqref{eq:otto_metric}) of $E(\rho)$, or 
	\item $I(\rho)$ is the rate of change of entropy along the heat flow $\dot\rho = \Delta\rho$.
\end{itemize}
In our context, $I(\rho)$ plays the role of negative potential energy in the Lagrangian
\begin{equation}\label{eq:regularized_lag}
	L(\rho,\dot\rho) = \frac{1}{2}\pair{\dot\rho,\dot\rho}_\rho + \varepsilon^2 I(\rho) .
\end{equation}
The corresponding action functional 
\begin{equation}\label{eq:action_regularized}
	A(\rho_t) = \int_0^1 \left( \frac{1}{2}\pair{\partial_t\rho_t,\partial_t\rho_t}_\rho + \varepsilon^2 I(\rho_t) \right)\ud t
\end{equation}
is called the \emph{entropic regularization} of \eqref{eq:action}.
Notice that the parameter $\varepsilon$ has physical dimension as thermal diffusivity.
This is not a coincidence.
Indeed, as we shall soon see this regularization significantly simplifies the variational problem (imaginary $\varepsilon$ also works!) by means of heat flows with thermal diffusivity $\varepsilon$.
Before that, however, let me just point out that the dynamical formulation of Benamou and Brenier~\cite{BeBr2000} is still applicable: if we change variables to $(\rho,m)$ as before we obtain again a convex optimization problem (since the functional $I(\rho)$ is convex)
\begin{equation}
	\min_{\rho_t,m_t} \int_0^1 \left( \frac{1}{2}\inner{m_t/\rho_t,m_t}_{L^2} + \varepsilon^2 I(\rho_t) \right)\ud t \quad\text{subject to}\quad \partial_t\rho_t + \divv(m_t) = 0 . 
\end{equation}
In a suitable setting one can apply convex analysis to obtain existence and uniqueness (\emph{cf.}~\cite{BeBr2000}).

So far I used Lagrangian mechanics to describe the dynamical formulation.
Let me now switch to the Hamiltonian view-point.
The Legendre transform of the Lagrangian \eqref{eq:regularized_lag} is
\begin{equation}\label{eq:legendre}
	T\Dens(M) \ni (\rho,\dot\rho) \mapsto \left(\rho, \frac{\delta L}{\delta \dot\rho}\right) = (\rho,S) \in T^*\Dens(M)
\end{equation}
where the co-tangent space $T^*_\rho\Dens(M)$ is given by co-sets of smooth functions defined up to addition of constants (if $M$ is connected as I assume). 
Notice that the $S$ in \eqref{eq:legendre} is exactly the $S$ in \eqref{eq:otto_metric} (defined only up to a constant).
The corresponding Hamiltonian on $T^*\Dens(M)$ is
\begin{equation}\label{eq:hamiltonian}
	H(\rho,S) = \frac{1}{2}\pair{S,S}_{*\rho} - \varepsilon^2 I(\rho),
\end{equation}
where the \emph{dual metric} $\pair{\cdot,\cdot}_*$ is given by
\begin{equation}
	\pair{S,S}_{*\rho} = \int_M \abs{\nabla S}^2\rho .
\end{equation}

Before I continue, consider a finite-dimensional analog of the Hamiltonian \eqref{eq:hamiltonian}: on $T^*\RR$ take $H(q,p) = p^2/2 - \varepsilon^2 q^2/2$. 
Of course, the analogy is $p^2/2 \leftrightarrow \pair{S,S}_{*\rho}$ and $q^2/2 \leftrightarrow I(\rho)$.
The equations of motion are
\begin{equation}
	\dot q = p, \quad \dot p = \varepsilon^2 q .
\end{equation}
This system describes a harmonic oscillator when $\varepsilon$ is \emph{imaginary}.
For real $\varepsilon$, the dynamics is not oscillatory.
Indeed, if we change variables to $x = p - \varepsilon q$ and $y = p + \varepsilon q$ we obtain the Hamiltonian $H(x,y) = xy/2$ with dynamics
\begin{equation}
	\dot x = x/2, \quad \dot y = -y/2.
\end{equation}
These are two uncoupled equations where $x$ is growing and $y$ is decaying exponentially.
Thus, we can expect this type of dynamics also for \eqref{eq:hamiltonian}.
Indeed, I shall now introduce a change of coordinates for $(\rho,S)$, analogous to the change of coordinates $(q,p)\iff (x,y)$.

\begin{definition}
	The \emph{imaginary Madelung transform}\footnote{In the standard Madelung transform $\varepsilon = \ii \hbar$ which is why I say `imaginary' here.}\ is given by 
	\begin{equation}
		T^*\Dens(M)\ni (\rho,S) \longmapsto \left( \underbrace{\sqrt{\rho\ee^{S/\epsilon}}}_{a}, \underbrace{\sqrt{\rho\ee^{-S/\epsilon}}}_{\bar b}  \right) ,
	\end{equation}
	where $a,\bar b \in C^\infty(M)$ are defined up to $(\ee^\sigma a , \ee^{-\sigma} \bar b)$ for $\sigma \in \RR$ and should fulfill $\int_M a \bar b = 1$.
	The individual component $a$ is known as the \emph{Hopf-Cole transform}.
\end{definition}

This transformation is a symplectomorphism (see \cite{Le2019} for $\varepsilon\in\RR$ and \cite{Re2012,KhMiMo2019} for $\varepsilon\in \ii \RR$).
The inverse transform is $\rho = a\bar b$ and $S = \varepsilon \log(a/\bar b)$.
The Hamiltonian~\eqref{eq:hamiltonian} expressed in the new canonical coordinates $(a,\bar b)$ thus become
\begin{equation}\label{eq:hamiltonian_after_madelung}
\begin{aligned}
	H(a,\bar b) &= \int_M \frac{\epsilon^2}{2} \left( \left| \frac{\nabla a}{a} - \frac{\nabla \bar b}{\bar b} \right|^2 - \left| \frac{\nabla a}{a} + \frac{\nabla \bar b}{\bar b} \right|^2\right) a\bar b	\\
	&= -\varepsilon^2 \int_M \nabla a \cdot \nabla \bar b = \varepsilon^2 \int_M  a \Delta \bar b \; .
\end{aligned}
\end{equation}
Notice two things: (i) the Hamiltonian is quadratic, and (ii) it is of the form in the toy example $H(x,y)$ above.
Hamilton's equations of motion for \eqref{eq:hamiltonian_after_madelung} are
\begin{equation}
	\dot a = \varepsilon\Delta a, \qquad \dot{\bar b} = -\varepsilon\Delta \bar b \, .
\end{equation}
Again, two decoupled equations as in the toy example, but now given by forward and backward heat flows with thermal diffusivity $\varepsilon$.

\begin{remark}\label{rem:non_uniqueness}
	To be more precise, one should also take into account that $(a,\bar b)$ is a co-set, so the general form of the equation should be
	\begin{equation}
		\dot a = \varepsilon\Delta a + \sigma a, \qquad \dot{\bar b} = -\varepsilon\Delta \bar b - \sigma \bar b \, ,
	\end{equation}
	where $t\mapsto \sigma(t)\in \RR$ is arbitrary. 
	However, since the scaling is arbitrary we can always represent the $(a,\bar b)$ co-set by the element for which $\sigma = 0$.
	Notice also that the constraint $\int_M a\bar b = 1$ is preserved by the flow, as a short calculation shows.
\end{remark}

It is, of course, not good to work with backward heat flows, but there is an easy fix.
Let $b(x,t) \coloneqq \bar b(x,1-t)$.
Then $b$ fulfills the forward heat equation (but backwards in time).
In the variables $a_t = a(\cdot,t)$ and $b_t = b(\cdot,t)$ the solution to the variational problem for the action \eqref{eq:action_regularized} must therefore fulfill
\begin{equation}\label{eq:EL_ab}
	\left\{
	\begin{aligned}
		&\partial_t a_t = \varepsilon\Delta a_t, & &a_0 b_1 = \rho_0\\
		&\partial_t b_t = \varepsilon\Delta b_t, & &a_1 b_0 = \rho_1 .
	\end{aligned}
	\right.
\end{equation}
The two equations are coupled only through mixed boundary conditions at $t=0$ and $t=1$.
With $a\coloneqq a_0$ and $b\coloneqq b_0$ these equations can be written in terms of the heat semigroup as
\begin{equation}\label{eq:stationary_eq}
	a\ee^{\varepsilon\Delta}b = \rho_0, \quad b\ee^{\varepsilon\Delta}a = \rho_1 .
\end{equation}
As you can see, a solution to \eqref{eq:stationary_eq} is a stationary point of the integral equations \eqref{eq:integral_equations_ab}.
Indeed, one should think of \eqref{eq:integral_equations_ab} as a gradient-type flow for solving the equations~\eqref{eq:stationary_eq}, as I shall now elaborate on.

If we take only the first part of the equations~\eqref{eq:stationary_eq} we obtain the equation 
\begin{equation}\label{eq:integral_equations_a}
		\partial_s a = -a \log\left(\frac{a\ee^{\varepsilon \Delta}b}{\rho_0} \right),
\end{equation}
with $b$ now as a fixed parameter.
Let $\sigma = a \ee^{\varepsilon\Delta}b$, so that $\partial_s \sigma = \ee^{\varepsilon\Delta}b \partial_s a$ (since $b$ is considered constant).
The Fisher-Rao metric for $\partial_s\sigma$ is given by 
\begin{equation}\label{eq:FRa}
	\langle \partial_s\sigma,\partial_s\sigma\rangle_{\sigma} = \int_M \frac{(\partial_s \sigma)^2}{\sigma} = \int_M \frac{(\partial_s a)^2\ee^{\varepsilon\Delta}b}{a}.
\end{equation}
Furthermore, the entropy of $\sigma$ relative to $\rho_0$ is given by 
\begin{equation*}
	H_{\rho_0}(\sigma) = \int_M \sigma \log\left(\frac{\sigma}{\rho_0} \right).
\end{equation*}
It is straightforward to check that the Riemannian gradient flow for the functional $F_1(a) = H_{\rho_0}(a\ee^{\varepsilon\Delta}b)$ with respect to the metric \eqref{eq:FRa} is given by equation~\eqref{eq:integral_equations_a}.
Likewise, the equation for $b$ with fixed $a$ is the Fisher-Rao gradient flow of $F_2(b) = H_{\rho_1}(b \ee^{\varepsilon\Delta}a)$.
Consequently, the Sinkhorn algorithm is the composition of steps for the first and second gradient flows.
The question of assigning one Riemannian gradient structure to the entire flow is more intricate, since the functionals $F_1$ and $F_2$ depend on both $a$ and $b$.

\appendix

\section{Beurling's ``forgotten'' result}

Motivated by Einstein's work on Brownian motion governed in law by the heat flow, Schrödinger~\cite{Sc1931} arrived at the equations \eqref{eq:stationary_eq} by studying the most likely stochastic path for a system of particles from initial distribution~$\rho_0$ to final distribution~$\rho_1$.
He gave physical arguments for why the problem should have a solution, but mathematically it was left open. 
S.~Bernstein then addressed it at the 1932 International Congress of Mathematics in Zürich. 
A full resolution, however, did not come until 1960 through the work of Beurling~\cite{Be1960}.
The objective was, in Beurling's own words, ``to derive general results concerning systems
like \eqref{eq:stationary_eq} and, in particular, to disclose the inherent and rather simple nature of the problem.''
Beurling certainly succeeded in doing so.
But to his astonishment (and slight annoyance) no-one took notice.
In fact, Schrödinger's bridge problem was itself largely forgotten among physicists and mathematicians.
Both Schrödinger's problem and the solution by Beurling were ``rediscovered'' and advocated by Zambrini~\cite{Za1986} as he was working with an alternative version of Nelson's framework for \emph{stochastic mechanics}~(cf.~\cite{Ne2011}).

Beurling relaxed the problem~\eqref{eq:stationary_eq} by replacing the functions $\rho_0,\rho_1,a,b$ by measures $\mu_0,\mu_1,\alpha,\beta$ on $M$. 
By multiplying the right-hand sides, one obtains the product measure $\mu \equiv \mu_0\wedge \mu_1$ on $M\times M$.
For the left-hand side, Beurling went on as follows. 
Any measure $\nu$ on $M\times M$ gives rise to the generalized marginal measures $\nu_0$ and $\nu_1$ defined for all $h\in C_0(M)$ by
\begin{equation}
	\int_M h \,d \nu_0 = \int_{M\times M} K_\epsilon(x,y) h(x) \, d\nu \quad\text{and}\quad
	\int_M h \,d \nu_1 = \int_{M\times M} K_\epsilon(x,y) h(y) \, d\nu .
\end{equation}
Thus, we have a mapping from the space of measures to the space of product measures via the quadratic map
\begin{equation}
	T_\epsilon\colon \nu \mapsto \nu_0 \wedge \nu_1 .
\end{equation}
Beurling noticed that the generalized version of Schrödinger's problem in equation~\eqref{eq:stationary_eq} can be written 
\begin{equation}\label{eq:schrodinger_beurling}
	T_\epsilon(\alpha\wedge\beta) = \mu_0\wedge\mu_1 . 
\end{equation}

Let $\mathcal{M}$ denote the space of Radon measures on $M\times M$
(i.e., the continuous dual of compactly supported continuous functions on $M\times M$) and $\mathcal{P} \subset \mathcal{M}$ the sub-set of product measures.
Further, let $\mathcal{M}^+\subset \mathcal{M}$ and $\mathcal{P}^+\subset\mathcal{P}$ denote the corresponding sub-sets of non-negative measures.
Since $K_\epsilon > 0$, it follows that
\begin{equation}\label{eq:beurling_map}
	T_\epsilon\colon \mathcal{M}^+ \to \mathcal{P}^+ .
\end{equation}
Let me now state Beurling's result adapted to the setting here.\footnote{The result proved by Beurling is much more general: it solves the problem for an $n$-fold product measure on the Cartesian product of $n$ locally compact Hausdorff spaces.}

\begin{theorem}[Beurling~\cite{Be1960}, Thm.~I]
	Let $M$ be compact (possibly with boundary) and $\epsilon > 0$.
	Then the mapping \eqref{eq:beurling_map} restricted to $\mathcal P^+$ is an automorphism (in the strong topology of $\mathcal M$).
\end{theorem}

From this result, a solution to Schrödinger's problem~\eqref{eq:schrodinger_beurling} in the category of measures is obtained as $\alpha\wedge\beta = T_\epsilon^{-1}(\mu_0\wedge\mu_1)$. Furthermore, the solution $\alpha\wedge\beta$ depends continuously (in operator norm) on the data $\mu_0\wedge\mu_1$.
Notice that, whereas $\alpha\wedge\beta$ is unique as a product measure, the components $\alpha,\beta$ themselves are only defined up to multiplication $e^f\alpha,e^{-f}\beta$ by an arbitrary function $f$ on $M$.
Thus, to work with product measures naturally captures the non-uniqueness pointed out in Remark~\ref{rem:non_uniqueness} above.  

The condition that $M$ is compact is used to obtain a positive lower and upper bound on the kernel $K_\epsilon$ (these are, in fact, the only conditions that Beurling's proof imposes on $K_\epsilon$). 
Such bounds are necessary for the map $T_\epsilon$ to be an automorphism (i.e., continuous with continuous inverse).
Beurling also gave a second, weaker result, which can be applied to the case of non-compact $M$.

\begin{theorem}[Beurling~\cite{Be1960}, Thm.~II]
	Let $\epsilon > 0$ and let $\mu_0\wedge \mu_1 \in \mathcal{P}^+$ be such that
	\begin{equation}
		\Big| \int_{M}\int_M \log K_\epsilon \, d\mu_0\, d\mu_1 \Big| < \infty .
	\end{equation}
	Then there exists a unique non-negative product measure $\nu$ on $M\times M$ that solves the equation 
	\begin{equation*}
		T_\epsilon(\nu) = \mu_0\wedge \mu_1 .
	\end{equation*}
\end{theorem}

Beurling's results can be viewed as a generalization from matrices to measures of Sinkhorn's theorem~\cite{Si1964} on doubly stochastic matrices, only it came four years \emph{before} Sinkhorn's result.
I find it remarkable that Beurling came up with these results independently of Kantorovich's formulation of optimal transport in terms of measures on a product space (which came to general knowledge in the West in the late 1960's).






\bibliographystyle{amsplainnat} 
\bibliography{/Users/moklas/Documents/Papers/References}

\end{document}

%% file: definitions.tex
\DeclareMathOperator{\divv}{div}

\newcommand{\ee}{\mathrm{e}}
\newcommand{\ii}{\mathrm{i}}
\newcommand{\pair}[1]{\left\langle #1 \right\rangle}

\newcommand{\inner}[1]{\langle\!\langle #1 \rangle\!\rangle}
\providecommand{\norm}[1]{\lVert#1\rVert}
\providecommand{\abs}[1]{\lvert#1\rvert}
\providecommand{\vect}[1]{\boldsymbol{#1}}

\newcommand{\ud}{\mathrm{d}}

\newcommand{\uud}{\,\ud}

\newcommand{\RR}{{\mathbb R}}

\newcommand{\Dens}{\mathrm{Dens}}







%% file: sinkhorn_geometry.bbl
\def\cprime{$'$}
\begin{thebibliography}{25}
\providecommand{\natexlab}[1]{#1}
\providecommand{\url}[1]{\texttt{#1}}
\providecommand{\urlprefix}{URL }
\providecommand{\eprint}[2][]{\url{#2}}

\bibitem[{Arnold and Khesin(1998)}]{ArKh1998}
V.~I. Arnold and B.~Khesin, \emph{Topological Methods in Hydrodynamics},
  Springer-Verlag, New York, 1998.

\bibitem[{Benamou and Brenier(2000)}]{BeBr2000}
J.-D. Benamou and Y.~Brenier, A computational fluid mechanics solution to the
  {M}onge--{K}antorovich mass transfer problem, \emph{Numer. Math.} \textbf{84}
  (2000), 375--393.

\bibitem[{Berman(2020)}]{Be2020}
R.~J. Berman, The {S}inkhorn algorithm, parabolic optimal transport and
  geometric {M}onge--{A}mp{\`e}re equations, \emph{Numer. Math.} \textbf{145}
  (2020), 771--836.

\bibitem[{Beurling(1960)}]{Be1960}
A.~Beurling, An automorphism of product measures, \emph{Ann. of Math.}
  \textbf{72} (1960), 189--200.

\bibitem[{Cuturi(2013)}]{Cu2013}
M.~Cuturi, Sinkhorn distances: Lightspeed computation of optimal transport,
  \emph{Advances in Neural Information Processing Systems 26}, pp. 2292--2300,
  Curran Associates, Inc., 2013.

\bibitem[{Feydy(2020)}]{Fe2020_thesis}
J.~Feydy, \emph{Analyse de donn{\'e}es g{\'e}om{\'e}triques, au del{\`a} des
  convolutions}, Ph.D. thesis, {Universit{\'e} Paris-Saclay}, 2020.

\bibitem[{Hamilton(1982)}]{Ha1982}
R.~S. Hamilton, The inverse function theorem of {N}ash and {M}oser, \emph{Bull.
  Amer. Math. Soc. (N.S.)} \textbf{7} (1982), 65--222.

\bibitem[{Khesin et~al.(2019)Khesin, Misio{\l}ek, and Modin}]{KhMiMo2019}
B.~Khesin, G.~Misio{\l}ek, and K.~Modin, Geometry of the {M}adelung transform,
  \emph{Arch. Rational Mech. Anal.} \textbf{234} (2019), 549--573.

\bibitem[{Khesin and Wendt(2009)}]{KhWe2009}
B.~Khesin and R.~Wendt, \emph{The Geometry of Infinite-dimensional Groups},
  vol.~51 of \emph{A Series of Modern Surveys in Mathematics}, Springer-Verlag,
  Berlin, 2009.

\bibitem[{L{\'e}ger(2019)}]{Le2019}
F.~L{\'e}ger, A geometric perspective on regularized optimal transport,
  \emph{Journal of Dynamics and Differential Equations} \textbf{31} (2019),
  1777--1791.

\bibitem[{L{\'e}ger and Li(2021)}]{LeLi2021}
F.~L{\'e}ger and W.~Li, {H}opf--{C}ole transformation via generalized
  {S}chr{\"o}dinger bridge problem, \emph{J. Differential Equations}
  \textbf{274} (2021), 788--827.

\bibitem[{Lehmann et~al.(2021)Lehmann, von Renesse, Sambale, and
  Uschmajew}]{LeReSaUs2021}
T.~Lehmann, M.-K. von Renesse, A.~Sambale, and A.~Uschmajew, A note on
  overrelaxation in the {S}inkhorn algorithm, \emph{Optimization Lett.}
  \textbf{16} (2021), 2209--2220.

\bibitem[{Leonard(2014)}]{Le2014}
C.~Leonard, A survey of the {S}chr{\"o}dinger problem and some of its
  connections with optimal transport, \emph{Discrete Contin. Dyn. Syst.}
  \textbf{34} (2014).

\bibitem[{McLachlan and Quispel(2002)}]{McQu2002}
R.~I. McLachlan and G.~R.~W. Quispel, Splitting methods, \emph{Acta Numer.}
  \textbf{11} (2002), 341--434.

\bibitem[{Modin(2017)}]{Mo2017}
K.~Modin, Geometry of matrix decompositions seen through optimal transport and
  information geometry, \emph{J. Geom. Mech.} \textbf{9} (2017), 335--390.

\bibitem[{Nelson(2012)}]{Ne2011}
E.~Nelson, Review of stochastic mechanics, \emph{Journal of Physics: Conference
  Series} \textbf{361} (2012), 012011.

\bibitem[{Otto(2001)}]{Ot2001}
F.~Otto, The geometry of dissipative evolution equations: the porous medium
  equation, \emph{Comm. Partial Differential Equations} \textbf{26} (2001),
  101--174.

\bibitem[{Peyr{\'e} et~al.(2019)Peyr{\'e}, Chizat, Vialard, and
  Solomon}]{PeChViSo2019}
G.~Peyr{\'e}, L.~Chizat, F.-X. Vialard, and J.~Solomon, Quantum entropic
  regularization of matrix-valued optimal transport, \emph{Eur. J. Appl. Math.}
  \textbf{30} (2019), 1079--1102.

\bibitem[{Peyr{\'e} and Cuturi(2020)}]{PeCu2020}
G.~Peyr{\'e} and M.~Cuturi, Computational optimal transport, \emph{Foundations
  and Trends in Machine Learning} \textbf{11} (2020), 355--607.

\bibitem[{Schmitzer(2019)}]{Sc2019}
B.~Schmitzer, Stabilized sparse scaling algorithms for entropy regularized
  transport problems, \emph{SIAM J. Sci. Comput.} \textbf{41} (2019),
  A1443--A1481.

\bibitem[{Schr{\"o}dinger(1931)}]{Sc1931}
E.~Schr{\"o}dinger, \emph{{\"U}ber die umkehrung der naturgesetze}, Verlag der
  Akademie der Wissenschaften in Kommission bei Walter De Gruyter, 1931.

\bibitem[{Sinkhorn(1964)}]{Si1964}
R.~Sinkhorn, A relationship between arbitrary positive matrices and doubly
  stochastic matrices, \emph{Ann. Math. Stat.} \textbf{35} (1964), 876--879.

\bibitem[{Thibault et~al.(2021)Thibault, Chizat, Dossal, and
  Papadakis}]{ThChDoPa2021}
A.~Thibault, L.~Chizat, C.~Dossal, and N.~Papadakis, Overrelaxed
  {S}inkhorn--{K}nopp algorithm for regularized optimal transport,
  \emph{Algorithms} \textbf{14} (2021), 143.

\bibitem[{von Renesse(2012)}]{Re2012}
M.-K. von Renesse, An optimal transport view of {S}chr{\"o}dinger's equation,
  \emph{Canad. Math. Bull} \textbf{55} (2012), 858--869.

\bibitem[{Zambrini(1986)}]{Za1986}
J.~C. Zambrini, Variational processes and stochastic versions of mechanics,
  \emph{J. Math. Phys.} \textbf{27} (1986), 2307--2330.

\end{thebibliography}
